\documentclass[11pt]{amsart}
\usepackage{amsfonts,latexsym,rawfonts,amsmath,amssymb,amsthm}
\usepackage[colorlinks,linkcolor=blue]{hyperref} 
\usepackage{verbatim}                    
\usepackage{color}
\usepackage{extarrows}   

\newtheorem{theorem}{Theorem}[section]
\newtheorem{lemma}[theorem]{Lemma}
\newtheorem{corollary}[theorem]{Corollary}

\theoremstyle{definition}
\newtheorem{definition}[theorem]{Definition}

\theoremstyle{remark}
\newtheorem{remark}[theorem]{Remark}

\numberwithin{equation}{section}

\allowdisplaybreaks

\begin{document}

\title[H$\ddot{\mbox{o}}$lder estimates and asymptotic behaviors]
{H$\ddot{\mbox{O}}$lder estimates and asymptotic behaviors at infinity of solutions of some degenerate elliptic equations in half spaces}

\author[X.B. Jia \and S.S. Ma]
{Xiaobiao Jia \and Shanshan Ma}
\address{School of Mathematics and Statistics, North China University of Water Resources and Electric Power,
           Zhengzhou 450046, China.}
\email{jiaxiaobiao@ncwu.edu.cn}

\address{School of Mathematics and Statistics, Zhengzhou University,
           Zhengzhou 450001, China.}
\email{mass1210@zzu.edu.cn}

\thanks{The second author was partly supported by Postdoctoral Research Foundation of Henan Province (No. 202001009).}




\begin{abstract}
In this paper we investigate the asymptotic behaviors at infinity of viscosity solutions of a kind of degenerate elliptic equations, which are linearized equations of a class of degenerate Monge-Amp\`ere equations.
Meanwhile, the H$\ddot{\mbox{o}}$lder estimates up to the boundary will be obtained by using the rescaling method, and a Liouville type result on some kind of Baouendi-Grushin type operator was deduced as a byproduct.
\end{abstract}

\maketitle

\section{Introduction}

In this paper we investigate the asymptotic behavior at infinity of viscosity solutions of the following degenerate non-divergence elliptic equations
\begin{equation}\label{main-Eq}
  Lu=x_n^{2\alpha}\sum_{i,j=1}^{n-1}a_{ij}(x)D_{ij}u(x)
       +2x_n^{\alpha}\sum_{i=1}^{n-1}a_{in}(x)D_{in}u(x)+D_{nn}u(x)=0\quad\mbox{in } \mathbb{R}^n_+\backslash \overline{B}_1^+,
\end{equation}
where $n\geq2$, $\alpha>0$, $\mathbb{R}^n_+=\mathbb{R}^n\cap\{x_n>0\}$, $B_1^+=\mathbb{R}^n_+\cap\{|x|<1\}$. Note that \eqref{main-Eq} relates to the linearization of the following degenerate Monge-Amp\`{e}re equation
\begin{equation}\label{EQ-monge-ampere}
    \det D^2u=f(x)x_n^{2\alpha}\quad \mbox{on }\{x_n>0\},
\end{equation}
where $\alpha>0$, and $f(x)$ is positive and continuous.

To ensure the ellipticity of operator $L$, it is always assumed that $a_{ij}(x), a_{in}(x)\in C(\mathbb{R}^n_+)$ ($i,j=1,\cdots,n-1$) and there exist constants $0<\lambda\leq\Lambda<\infty$
such that for any $\xi\in \mathbb{R}^{n-1}$,
\begin{equation}\label{Sz-unif-ellp}
  \lambda|\xi|^2\leq {\xi}^T\sum_{i,j=1}^{n-1}a_{ij}(x)\xi \leq \Lambda |\xi|^2, \quad\forall\; x\in\mathbb{R}^n_+,
\end{equation}
and for some $0<\delta<1$,
\begin{equation}\label{Sz-unif-ellp-2}
  1-\lambda^{-1}\sum_{i=1}^{n-1}||a_{in}||^2_{L^{\infty}(\mathbb{R}^n_+)}>\delta.
\end{equation}

In this paper, solutions always indicate viscosity solutions (cf.\cite{Caffarelli-Cabre-fully-nonlinear,Crandall-Ishii-Lions-1992-user,Han-Lin-EPDE-2011,
Ramaswamy-Ramaswamy-1996-Maximum} for definition).

 As for $\alpha=0$, by \eqref{Sz-unif-ellp} and \eqref{Sz-unif-ellp-2}, $L$ is uniformly elliptic. The asymptotic behavior at infinity was considered in \cite{Jia-Li-Li-2018-AdM}.
 Note that the crucial key to obtain the asymptotic behavior is the boundary H$\ddot{\mbox{o}}$lder estimates,
 which is classical for uniformly elliptic equations (cf.\cite{Caffarelli-Cabre-fully-nonlinear,Gilbarg-Trudinger-2001-Springer,Han-Lin-EPDE-2011}).

 As for $\alpha>0$, When $a_{ij}\equiv 1$ and $a_{in}\equiv0$ ($i$, $j\leq n-1$), $L$ is a
Baouendi-Grushin type operator,  which was introduced by \cite{Baouendi-1966} and \cite{Grushin-1970}, respectively.  There have been extensive works on the studies of the Baouendi-Grushin type operators (cf.\cite{Bauer-Furutani-Iwasaki-2015-AdV,Franchi-Lanconelli-1983,Garofalo-1993-JDE,Kombe-2006-MN, Monti-Morbidelli-2006-Duke,Robinson-Sikora-2008-MZ} and references therein). When $a_{ij}$ satisfies \eqref{Sz-unif-ellp}, Le and Savin \cite{Le-Savin-2017-Invn} obtained the boundary Schauder estimates for solutions of the degenerate elliptic equation

\begin{equation*}
  x_n^{\alpha}\sum_{i,j=1}^{n-1}a_{ij}(x)D_{ij}u(x)
  +D_{nn}u(x)=x_n^{\alpha}f(x)
  \quad\mbox{in } B_1^+
\end{equation*}
with $\alpha>0$.

 In this paper, we would like to study  the asymptotic behavior at infinity of solutions of \eqref{main-Eq} with the coefficients satisfying \eqref{Sz-unif-ellp} and \eqref{Sz-unif-ellp-2}.
The asymptotic behavior at infinity of solutions of \eqref{main-Eq} is a key point to obtain the asymptotic behavior at infinity of solutions of \eqref{EQ-monge-ampere} provided that $f(x)$ tends to some positive constant at infinity with proper decay rate.

By the similar rescaling method in \cite{Le-Savin-2017-Invn},  we establish the H$\ddot{\mbox{o}}$lder estimates up to the flat boundary of solutions of \eqref{main-Eq} with the coefficients satisfying \eqref{Sz-unif-ellp} and \eqref{Sz-unif-ellp-2}.
\begin{theorem}\label{thm1}
  Let $u\in C(\overline{B}_1^+)$ be a solution of
\begin{eqnarray}\label{Eq-holder}
       \left\{
        \begin{aligned}
     &Lu(x)=0\quad\mbox{in } B_1^+\\
    &u(x)=0\quad \mbox{ on }B_1\cap\{x_n=0\},
        \end{aligned}
     \right.
\end{eqnarray}
where $L$ is given by \eqref{main-Eq} with the coefficients satisfying \eqref{Sz-unif-ellp} and \eqref{Sz-unif-ellp-2}. Then $u\in C^{\frac{1}{1+\alpha}}(\overline{B}_{\frac{1}{2}}^+)$.
\end{theorem}

As we all known (cf. \cite{Garofalo-1993-JDE}), {the function}
\begin{equation}\label{def-dx}
    d(x)=d(x',x_n)=\left\{|x'|^2+\frac{1}{(\alpha+1)^2}x_n^{2(\alpha+1)}\right\}^{\frac{1}{2(\alpha+1)}}
\end{equation}
is the nature gauge associated with the Baouendi-Grushin type operator
\begin{equation}\label{oper-good}
    \mathfrak{L}=x_n^{2\alpha}\Delta_{x'}+\partial_{x_nx_n}^2\quad\mbox{in }\mathbb{R}^n_+.
\end{equation}
That is,
\[
   \mathfrak{ L}(d(x)^Q)=0\quad\mbox{in }\mathbb{R}^n_+,
\]
where
$
Q=(\alpha+1)(n-1)+1,
$
and see \emph{Section \ref{S3}} below for details.

\bigskip

Theorem \ref{thm1} together with Harnack inequalities and the comparison principle yields our main theorem as {follows}.
\begin{theorem}\label{Main-thm}
Let $u\in C^{1}(\mathbb{\overline{R}}^n_+\backslash B_1^+)$ be a solution of
\begin{equation}\label{EQ_Linear_Eq}
     \left\{
        \begin{aligned}
              &Lu=0\quad\mbox{in } \mathbb{R}^n_+\backslash \overline{B}_1^+, \\
              &u=0\quad\mbox{ on } \{x_n=0,|x|\geq 1\},\\
          \end{aligned}
     \right.
\end{equation}
where  $L$ is given by \eqref{main-Eq} with the coefficients satisfying \eqref{Sz-unif-ellp} and  \eqref{Sz-unif-ellp-2}, and for some $s>0$,
\begin{equation}\label{Sz-aij-good}
  |a_{ij}(x)-\delta_{ij}|+|a_{in}(x)|\leq  d(x)^{-s}\quad \mbox{in }\mathbb{R}^n_+\backslash \overline{B}_{1}^+,\quad i,j<n.
\end{equation}
Assume that $|u|\leq1$ on $\partial B_1\cap\{x_n>0\}$,
$|Du|\leq 1$ in $\mathbb{\overline{R}}^n_+\backslash B_1^+$
and $|Du|\rightarrow 0$ as $|x|\rightarrow \infty$.
Then
\begin{equation}\label{SZ_AsB_Lin}
    |u(x)|\leq \frac{C x_n}{d(x)^Q}
    \quad\mbox{in }    \mathbb{\overline{R}}^n_+\backslash B_R^+,
\end{equation}
where $d(x)$ is as in \eqref{def-dx}, $C>0$ and $R\geq1$ depend only on $\alpha$, $\delta$, $s$ and $n$.
\end{theorem}

\begin{remark}
If $\alpha=0$, Theorem \ref{Main-thm} still holds (cf.\cite{Jia-Li-Li-2018-AdM}).
\end{remark}

By Theorem \ref{Main-thm} and comparison principle, we have the following Liouville type theorem.

\begin{theorem}\label{tm-Liouville}
Let $u\in C^{1}(\mathbb{ {R}}^n_+ )$ be a solution of
\begin{equation}\label{EQ_Linear_Eq}
     \left\{
        \begin{aligned}
              &\mathfrak{L}u=0\quad\mbox{in } \mathbb{R}^n_+, \\
              &u=0\quad\mbox{ on } \{x_n=0\},\\
          \end{aligned}
     \right.
\end{equation}
where $\mathfrak{L}$ was given by \eqref{oper-good}.
If $|Du|\rightarrow 0$ as $|x|\rightarrow \infty$.
Then $u(x)$ must be zero.
\end{theorem}

Note that the proof of Theorem \ref{tm-Liouville} is obviously, and thus we omit it here.

\bigskip

This paper is organized as follows. In \emph{Section \ref{S2}}, we first claim that $L$ is uniformly elliptic in interior, and then the interior H$\ddot{\mbox{o}}$lder estimates is clear.  {Secondly}, we show the boundary H$\ddot{\mbox{o}}$lder estimates,  {which} can be approached by the interior H$\ddot{\mbox{o}}$lder estimates via rescaling. In \emph{Section \ref{S3}}, a supersolution is constructed according to the fundamental solution of one Baouendi-Grushin type operator in the upper half space.  {Then it together with the H$\ddot{\mbox{o}}$lder estimates up to the flat boundary implies that Theorem \ref{Main-thm} holds.}

\section{Proof of Theorem \ref{thm1}}\label{S2}

Firstly, by \eqref{Sz-unif-ellp} and \eqref{Sz-unif-ellp-2},
we show that $L$ is elliptic in $\overline{B}_{1}^+$ and uniformly elliptic  {if $x_n>\varepsilon_0$ with any fixed $\varepsilon_0>0$}. Precisely,
\begin{lemma}\label{Lm-Unif-Ellp}
Let the coefficients of $L$ in \eqref{main-Eq} satisfy \eqref{Sz-unif-ellp} and \eqref{Sz-unif-ellp-2}. Then $L$ is elliptic in $\overline{B}_{1}^+$. Furthermore, for any fixed $\varepsilon_0>0$, $L$ is uniformly elliptic in $\overline{B}_{1}^+\cap\{x_n\geq \varepsilon_0\}$.
\end{lemma}

\begin{remark}
   Lemma \ref{Lm-Unif-Ellp} is standard, and  {one can} see its proof in Appendix below.
\end{remark}

To show the H$\ddot{\mbox{o}}$lder estimates up to the boundary, we {need to give} some notions as follows (cf.\cite{Le-Savin-2017-Invn}).

\begin{definition}
  We define a distance $d_\alpha$ between point $y$ and point $z$ by
  \[
  d_\alpha(y,z):=|y'-z'|+\left|y_n^{1+\alpha}-z_n^{1+\alpha}\right|.
  \]
\end{definition}
Observe that the relation between $d_\alpha$ and the Euclidean distance satisfies:
\begin{equation}\label{sz-rela-da-d}
  c|y-z|^{1+\alpha}\leq d_\alpha(y,z)\leq C|y-z|
\end{equation}
and
\begin{equation}\label{sz-rela-da-d-2}
 d_\alpha(y,z)\sim|y-z|\quad \mbox{if }y,z\in \overline{B}_{1}^+\cap\left\{x_n\geq \frac{1}{8}\right\}.
\end{equation}

  For any $h>0$ and any $\widetilde{x}\in \mathbb{R}^n$, we denote
\begin{equation}\label{Def-Eh}
  E_h(\widetilde{x})=\left\{x\in \mathbb{R}^n:|x'-\widetilde{x}'|^2+|x_n-\widetilde{x}_n|^{2(1+\alpha)}<h\right\},
\end{equation}
and
\[
F_h=\mbox{diag}\left(h^{\frac{1}{2}},h^{\frac{1}{2}},\cdots,h^{\frac{1}{2}},h^{\frac{1}{2(1+\alpha)}}\right).
\]
For simplicity, we always denote
\[
  E_h=E_h(0)=\left\{x\in \mathbb{R}^n:|x'|^2+|x_n|^{2(1+\alpha)}<h\right\};\quad  E_h^+= E_h\cap\{x_n>0\}.
\]

 {A} simple calculation gives
\begin{equation}\label{Sz-fheh}
F_h E_{\alpha'}\left(\frac{1}{2}e_n\right)
     =E_{\alpha'h}\left(\frac{1}{2}h^{\frac{1}{2(1+\alpha)}}e_n\right), \quad F_hE_1^+=E_h^+,
\end{equation}
where
\[
e_n=(0,\cdots,0,1), \quad\alpha'=4^{-2(1+\alpha)}.
\]

Note that \eqref{main-Eq} and $d_\alpha$ keep their forms under the transformation $x\rightarrow F_hx$. Precisely, {let}
\begin{equation}\label{Sz-transFh}
  \widetilde{u}(x)=u(F_hx),\quad x\in E_1,
\end{equation}
 { and then it} solves
\begin{equation}\label{eq-u-trans}
  \widetilde{L}\widetilde{u}=x_n^{2\alpha}\sum_{i,j=1}^{n-1}\widetilde{a}_{ij}(x)D_{ij}\widetilde{u}(x)
  +x_n^{\alpha}\sum_{i=1}^{n-1}2\widetilde{a}_{in}(x)D_{in}\widetilde{u}(x)+D_{nn}\widetilde{u}(x)=0
\end{equation}
with
\begin{equation}\label{Sz-coeff}
  \widetilde{a}_{ij}(x)=a_{ij}(F_hx),\quad \widetilde{a}_{in}(x)=a_{in}(F_hx),\quad i,j\leq n-1,
\end{equation}
and
\begin{equation}\label{sz-def-dada}
  d_{\alpha}(y,z)=h^{-\frac{1}{2}}d_{\alpha}(F_hy,F_hz).
\end{equation}

If function $w$ is $\gamma$-H$\ddot{\mbox{o}}$lder continuous in $\Omega\subset\overline{B}_1^+$ with respect to $d_\alpha$, we write
\[
w\in C^{\gamma}_\alpha(\overline{\Omega})
\]
and define
\[
[w]_{C^{\gamma}_\alpha(\overline{\Omega})}=\sup_{y,z\in\overline{\Omega},y\neq z}\frac{|w(y)-w(z)|}{(d_\alpha(y,z))^\gamma},\quad
||w||_{C^{\gamma}_\alpha(\overline{\Omega})}
=||w||_{L^{\infty}(\overline{\Omega})}
+[w]_{C^{\gamma}_\alpha(\overline{\Omega})}.
\]

\bigskip

Next we show Theorem \ref{thm1}.

  \emph{Proof of Theorem \ref{thm1}.}

We divided this into two case:
$u\in C^{\frac{1}{1+\alpha}}(\overline{B}_{\frac{1}{2}}^+\cap\{x_n>\frac{1}{8}\})$ and $u\in C^{\frac{1}{1+\alpha}}(\overline{B}_{\frac{1}{2}}^+\cap\{x_n\leq\frac{1}{8}\})$.

\emph{\textbf{Case 1.}} $x\in \overline{B}_{\frac{1}{2}}^+\cap\{x_n>\frac{1}{8}\}$.

By Lemma \ref{Lm-Unif-Ellp}, $L$ is uniformly elliptic in
$ \overline{B}_{\frac{1}{2}}^+\cap\{x_n>\frac{1}{8}\}$. Applying the classical H$\ddot{\mbox{o}}$lder estimates to $u$,
 there exists $C>0$, depending only on $\lambda$, $\Lambda$, $\alpha$, $\delta$, $n$ and $||u||_{L^{\infty}}$, such that
\[
     [u]_{C^{\frac{1}{1+\alpha}}\left(\overline{E_{\alpha'}\left(\frac{1}{2}e_n\right)}\right)}
    \leq C||u||_{L^{\infty}}\leq C.
\]

\emph{\textbf{Case 2.}} $x\in \overline{B}_{\frac{1}{2}}^+\cap\{x_n\leq\frac{1}{8}\}$.
We show this case by four steps as the following.

\emph{Step 1.}  {There} exists some $C>0$, depending only on $\lambda$, $\Lambda$, $\alpha$, $\delta$, $n$ and $||u||_{L^{\infty}}$, such that
\begin{equation}\label{ZS-u-xn}
|u(x)|\leq Cx_n \quad \mbox{in } B_{\frac{3}{4}}^+.
\end{equation}

It only  {need} to show that for any $x_0\in\{x_n=0,|x'|<\frac{3}{4}\}$,
\[
|u(x_0,x_n)|\leq Cx_n.
\]

Let
\[
\overline{u}(x)=Cx_n+B|x'-x_0'|^2-\frac{C}{2}x_n^{2+\alpha}
\]
with $B=16||u||_{L^{\infty}}$.
One can choose $C>0$, depending only on $\Lambda$, $\alpha$, $n$ and $||u||_{L^{\infty}}$, such that
\begin{equation}\label{EQ-comp-1}
L\overline{u}\leq 0\quad \mbox{in } B_1^+,\qquad
\overline{u}\geq ||u||_{L^{\infty}}\geq u\quad \mbox{on } \partial B_1^+\textcolor{red}{,}
\end{equation}
if taking
\[
2(n-1)\Lambda B-(2+\alpha)(1+\alpha)C/2\leq 0
\quad \mbox{and} \quad
\frac{C}{2}x_n+B|x'-x_0'|^2>||u||_{L^{\infty}}\quad \mbox{on } \partial B_1^+.
\]

Therefore,  \eqref{EQ-comp-1} and the comparison principe (cf.\cite[Theorem 6]{Ramaswamy-Ramaswamy-1996-Maximum}) yield \emph{Step 1}.

\emph{Step 2.}
For any fixed $h\in(0,1]$,
\begin{equation}\label{Sz-halfnorm-est}
   [u]_{
   C^{\frac{1}{1+\alpha}}_\alpha
   \left(\overline{E_{\alpha'h}
     \left(\frac{1}{2}h^{\frac{1}{2(1+\alpha)}}e_n\right)
     }\right)}\leq C.
\end{equation}

In fact, let $\widetilde{u}$ be as in \eqref{Sz-transFh},
and then $\widetilde{u}$ solves \eqref{eq-u-trans} in $B_1^+$.
By \eqref{Sz-transFh} and \eqref{ZS-u-xn},
we have
\begin{equation}\label{sz-hatu-h}
    \widetilde{u}\leq Ch^{\frac{1}{2(1+\alpha)}}\quad \mbox{in } B_1^+.
\end{equation}

Similar to \emph{\textbf{Case 1}}, applying the classical H$\ddot{\mbox{o}}$lder estimates to $\widetilde{u}$ in $E_{\frac{1}{4}}\left(\frac{1}{2}e_n\right)$, we have
\[
[\widetilde u]_{C^{\frac{1}{1+\alpha}}\left(\overline{E_{\alpha'}\left(\frac{1}{2}e_n\right)}\right)}
\leq C||\widetilde u||_{L^\infty(B_1^+)} \leq Ch^{\frac{1}{2(1+\alpha)}}.
\]
 {By \eqref{sz-rela-da-d-2}, we see}
 \[
 [\widetilde u]_{C^{\frac{1}{1+\alpha}}_\alpha\left(\overline{E_{\alpha'}\left(\frac{1}{2}e_n\right)}\right)}
  \leq Ch^{\frac{1}{2(1+\alpha)}}.
 \]
This together with \eqref{Sz-fheh}, \eqref{Sz-transFh} and \eqref{sz-def-dada} yields \eqref{Sz-halfnorm-est},
since
\[
\frac{|\widetilde{u}(y)-\widetilde{u}(z)|}{(d_\alpha(y,z))^{\frac{1}{1+\alpha}}}
=\frac{| u(F_hy)- u(F_hz)|}
{h^{-\frac{1}{2(1+\alpha)}}(d_\alpha(F_hy,F_hz))^{\frac{1}{1+\alpha}}}.
\]

\emph{Step 3.}
We prove that $u\in C^{\frac{1}{1+\alpha}}_{\alpha}$ at $0$ along $e_n$ direction, i.e.,
\[
\sup_{0<h<1}
\frac{\left|u\left(\frac{1}{2}h^{\frac{1}{2(1+\alpha)}}e_n\right)-u(0)\right|}
{\left((\frac{1}{2}h^{\frac{1}{2(1+\alpha)}})^{1+\alpha}\right)^{\frac{1}{1+\alpha}}}\leq C
\]
for some $C>0$ depending only on $\lambda$, $\Lambda$, $\alpha$ and $n$.

It suffices to prove that
\[
\left|u\left(\frac{1}{2}h^{\frac{1}{2(1+\alpha)}}e_n\right)-u(0)\right|
\leq C h^{\frac{1}{2(1+\alpha)}},
\]
where $C>0$ independents on $h$.

Indeed, \emph{Step 2} yields that for any $k=1,2,\cdots,$
\[
  \left|u\left(\frac{1}{2^k}h^{\frac{1}{2(1+\alpha)}}e_n\right)    -u\left(\frac{1}{2^{k+1}}h^{\frac{1}{2(1+\alpha)}}e_n\right)\right|
  \leq C2^{-k-1}h^{\frac{1}{2(1+\alpha)}},
\]
These imply that

\begin{align*}
     \frac{\left|u\left(\frac{1}{2}h^{\frac{1}{2(1+\alpha)}}e_n\right)-u(0)\right|} {h^{\frac{1}{2(1+\alpha)}}}
     \leq     \sum_{k=1}^{\infty}
     \frac{\left|u\left(\frac{1}{2^k}h^{\frac{1}{2(1+\alpha)}}e_n\right)
       -u\left(\frac{1}{2^{k+1}}h^{\frac{1}{2(1+\alpha)}}e_n\right)\right|}
       {h^{\frac{1}{2(1+\alpha)}}}
   \leq  \sum_{k=1}^{\infty}C2^{-k-1}\leq C.
  \end{align*}
Therefore, $u\in C^{\frac{1}{1+\alpha}}_{\alpha}$ at $0$ along $e_n$ direction.

\emph{Step 4.}
We show \textbf{\emph{Case 2}}.

Similar to \emph{Step 3}, we have that $u\in C^{\frac{1}{1+\alpha}}_{\alpha}$ at any $x\in B_{\frac{1}{2}}^+\cap\{x_n=0\}$ along $e_n$ direction.

Let $y,z\in \overline{B}_{\frac{1}{2}}^+\cap\{x_n\leq\frac{1}{8}\}$ and denote by $y_n$, $z_n$ the $n^{th}$ component of $y$ and $z$, respectively.

If $z\in E_{2^{-2(1+\alpha)}y_n^{2(1+\alpha)}}\left(y_n\right)$ or $y\in E_{2^{-2(1+\alpha)}z_n^{2(1+\alpha)}}\left(z_n\right)$,  by \eqref{Sz-halfnorm-est} in \emph{Step 2}, we are done.

Otherwise, $z\notin E_{2^{-2(1+\alpha)}y_n^{2(1+\alpha)}}\left(y_n\right)$ and $y\notin E_{2^{-2(1+\alpha)}z_n^{2(1+\alpha)}}\left(z_n\right)$, and then we have
\begin{equation}\label{sz-y-z-large}
  |y-z|^2\geq \max \left\{2^{-2(1+\alpha)}z_n^{2(1+\alpha)}, 2^{-2(1+\alpha)}y_n^{2(1+\alpha)}\right\}.
\end{equation}

By \emph{Step 3} and the boundary value condition, we get

\begin{align}\label{sz-1}
  |u(y)-u(z)|&\leq |u(y)-u(y',0)|+|u(y',0)-u(z',0)|+|u(z',0)-u(z)|\\
             &\leq C|y_n|+C|z_n|
             \leq C|x-y|^{\frac{1}{1+\alpha}}\quad (\mbox{by \eqref{sz-y-z-large}}).\nonumber
\end{align}
It follows that $u\in C^{\frac{1}{1+\alpha}}\left(\overline{B}_{\frac{1}{2}}^+\cap\{x_n\leq\frac{1}{8}\}\right)$.

Therefore, \textbf{\emph{Case 1}} and \textbf{\emph{Case 2}} finish the proof of Theorem \ref{thm1}.
$\hfill\Box$

\section{Proof of Theorem \ref{Main-thm}}\label{S3}

In this section we divide the proof of Theorem \ref{Main-thm} into two steps as the following. In fact, \emph{Subsection \ref{S3-1}} gives the convergence at infinity of the solutions in Theorem \ref{Main-thm}, and then \emph{Subsection \ref{S3-2}} shows its asymptotic behavior at infinity with decay rate.

Recall that the symbols $F_h$, $E_h$ and $E_h^+$ are defined in \emph{Section \ref{S2}}.

\subsection{The convergence at infinity}\label{S3-1}
In the subsection we apply the H$\ddot{\mbox{o}}$lder estimates up to the flat boundary to show that the solution in Theorem \ref{Main-thm} converges at infinity.

Hereinafter, we say some constant is universal if it depends only on $\lambda$, $\Lambda$, $\alpha$, $\delta$ and $n$. {The} universal constant may change from line to line if necessary.

A straightforward corollary of the boundary H$\ddot{\mbox{o}}$lder estimates is the following.

\begin{corollary}\label{Lm_Int_Str_Small}
  Let $u\in C(\overline{E_{4R}^+\backslash E_R^+})$ be a solution of
\begin{equation}\label{EQ_supsolu}
     \left\{
        \begin{aligned}
           &Lu=0\quad\mbox{in } E_{4R}^+\backslash \overline{E}_R^+,\\
           &u\leq1 \quad\,\mbox{ on }\partial( E_{4R}^+\backslash \overline{E}_R^+)\cap\{x_n>0\},\\
           &u \leq\frac{1}{2}\quad\;\mbox{on }
              \partial(E_{4R}^+\backslash\overline{E}_R^+)\cap\{x_n=0\},\\
        \end{aligned}
     \right.
\end{equation}
where  $L$ is given by \eqref{main-Eq} with coefficients satisfying \eqref{Sz-unif-ellp} and \eqref{Sz-unif-ellp-2} in $E_{4R}^+\backslash \overline{E}_R^+$ for some $R>0$. Then there exists universal constant $c_0>0$ such that
\[
    u(x)\leq 1-c_0\quad \mbox{on } \partial{E}_{2R}\cap\{x_n\geq0\}.
\]
\end{corollary}

\begin{proof}
It only {need} to prove this corollary with $u(x) =\frac{1}{2}$ on $\partial(E_{4R}^+\backslash\overline{E}_R^+)\cap\{x_n=0\}$.
Otherwise, one can consider a supersolution
$v$ with $v(x) =\frac{1}{2}$ on $\partial(E_{4R}^+\backslash\overline{E}_R^+)\cap\{x_n=0\}$, and if it holds for $v$, by the comparison principle, so does for $u$.

  Let
  \[
        \widehat{u}(x)=u(F_{R}x),\quad x\in E_{4}^+\backslash \overline{E}_1^+.
  \]
  By the definitions of $F_R$ and $E_R^+$ in \emph{Section \ref{S2}}, we have
$
  F_{R}(E_{4}^+\backslash \overline{E}_1^+)=E_{4R}^+\backslash \overline{E}_R^+.
$
Then
\begin{equation}\label{EQ_supsolu-2}
     \left\{
        \begin{aligned}
           &\widetilde{L}\widehat{u}=0\quad\mbox{ in } E_{4}^+\backslash \overline{E}_1^+,\\
           &\widehat{u}\leq1 \quad\;\;\mbox{ on }\partial( E_{4}^+\backslash \overline{E}_1^+)\cap\{x_n>0\},\\
           &\widehat{u} =\frac{1}{2}\quad\;\mbox{ on }
              \partial(E_{4}^+\backslash\overline{E}_1^+)\cap\{x_n=0\},\\
        \end{aligned}
     \right.
\end{equation}
where $\widetilde{L}$ is given by \eqref{eq-u-trans}. Clearly, coefficients of $\widetilde{L}$ also satisfy \eqref{Sz-unif-ellp} and \eqref{Sz-unif-ellp-2} in $E_{4}^+\backslash \overline{E}_1^+$. Then by the third {equality} in \eqref{EQ_supsolu-2} and Theorem \ref{thm1}, there exists {a universal } $0<\tau\leq1$ such that
\begin{equation}\label{SZ-u-lager-1}
    \widehat{u}(x)\leq \frac{2}{3}
    \quad \mbox{on }\partial E_2\cap\{0\leq x_n\leq\tau\}.
\end{equation}

By the comparison principle,
we have $\widehat{u}\leq1$ in $E_{4}^+\backslash \overline{E}_1^+$.  $1-\widehat{u}$ satisfies
\[
\widetilde{L}(1-\widehat{u})=0\quad\mbox{ in } E_{4}^+\backslash \overline{E}_1^+.
\]
 Then by the interior Harnack inequality to $1-\widehat{u}$, there exists  {a universal } $C\geq1$ such that
\[
    C\inf\limits_{\partial E_2\cap\{x_n\geq\tau\}} (1-\widehat{u})
    \geq\sup\limits_{\partial E_2\cap\{x_n\geq\tau\}} (1-\widehat{u})
    \geq\sup\limits_{\partial E_2\cap\{x_n=\tau\}} (1-\widehat{u})
    \geq \frac{1}{3}.
\]
This implies
\begin{equation}\label{SZ-u-lager-2}
   \widehat{u}(x)\leq 1-\frac{1}{3C}
   \quad \mbox{on }\partial E_2\cap\{x_n\geq\tau\}.
\end{equation}

This corollary follows from the definition of $\widehat{u}$, \eqref{SZ-u-lager-1} and \eqref{SZ-u-lager-2}, if we take  $c_0=\frac{1}{3C}$.
\end{proof}

Applying Corollary \ref{Lm_Int_Str_Small}, we have a convergence result as follows.

\begin{theorem}\label{Tm-converges}
  Let $u\in  C(\mathbb{\overline{R}}^n_+\backslash E_1^+)$ be a solution of
  \begin{equation*}
      Lu=0\quad\mbox{in } \mathbb{R}^n_+\backslash \overline{E}_1^+,
\end{equation*}
where  $L$ is given by \eqref{main-Eq} with the coefficients satisfying \eqref{Sz-unif-ellp} and \eqref{Sz-unif-ellp-2} in $\mathbb{R}^n_+\backslash \overline{E}_1^+$. If
\begin{itemize}
  \item  $|u|\leq1$  on $(\partial E_1\cap\{x_n>0\})\cup\{x_n=0,|x|\geq 1\}$,
  \item  $u(x',0)\rightarrow \beta$ as $|x'|\rightarrow\infty$
  \item  $|Du(x)|\rightarrow 0$ as $|x|\rightarrow \infty$.
\end{itemize}
Then $u(x)\rightarrow \beta$ as $|x|\rightarrow\infty$.
\end{theorem}
\begin{proof}
  The proof of this theorem is divided into two steps as follows.

\textbf{Step 1.} $|u|\leq1$ in $\mathbb{\overline{R}}^n_+\backslash \overline{E}_1^+$.

For any $\varepsilon>0$, {since} $|Du|\rightarrow 0$ as $|x|\rightarrow\infty$,
there exists $R_\varepsilon\geq1$ such that
\begin{equation}\label{SZ-Du-small}
  |Du|\leq \varepsilon
  \quad \mbox{in } \mathbb{\overline{R}}^n_+\backslash Q_{R_\varepsilon}^+,
\end{equation}
where
$Q_{R_\varepsilon}^+:=\{(x',x_n):|x'|<R_\varepsilon, 0<x_n<R_\varepsilon\}$
is a cylinder.

By $|u|\leq 1$ on $\{x_n=0,~|x|\geq 1\}$,
\eqref{SZ-Du-small} and Newton-Leibniz formula, we have
\[
    |u(x)|\leq 1+2\varepsilon x_n
    \quad \mbox{on } \partial Q_{R_\varepsilon}^+\cap\{x_n>0\}.
\]
Since $|u|\leq1$ on $(\partial E_{1}\cap\{x_n>0\})\cup\{x_n=0,|x|\geq 1\}$, we get
\[
    |u(x)|\leq 1+2\varepsilon x_n
    \quad\mbox{on } \partial(Q_{R_\varepsilon}^+\backslash \overline{E}_{1}^+).
\]

Obviously, $1+2\varepsilon x_n$ solves  \eqref{main-Eq}
in $Q_{R_\varepsilon}^+\backslash \overline{E}_{1}^+$. Then by the comparison principle, we get
\[
    |u(x)|\leq 1+2\varepsilon x_n
    \quad\mbox{in } Q_{R_\varepsilon}^+\backslash \overline{E}_{1}^+.
\]
Letting $\varepsilon\rightarrow 0$,
it completes the proof of this step.

\textbf{Step 2.} $u(x)\rightarrow \beta$ as $|x|\rightarrow\infty$.

Without loss of generality, we suppose that $\beta=0$. Otherwise, we consider
$
    \frac{u(x)-\beta}{1+|\beta|}.
$

Now we argue by contradiction.
If this step is not true, by Step 1,
$u$ has finite superior limit $\overline{u}>0$ or inferior limit $\underline{u}<0$ at infinity.
It suffices to assume $\overline{u}>0$.

By the definition of $\overline{u}$ and $u(x',0)\rightarrow \beta$ as $|x'|\rightarrow\infty$,
there exists large $R_1\geq1$ such that for all $R\geq R_1$,
\[
    u(x)\leq\left(1+\frac{c_0}{2}\right)\overline{u}
    \quad\mbox{in }\mathbb{\overline{R}}^n_+\backslash \overline{E}_R^+
\]
and
\[
    u(x',0)\leq \frac{1}{2}\left(1+\frac{c_0}{2}\right)\overline{u}
    \quad\mbox{if }|x'|\geq R,
\]
where $c_0$ is given by Corollary \ref{Lm_Int_Str_Small}.
Then applying Corollary \ref{Lm_Int_Str_Small}
to
$
\frac{u(x)}{\left(1+\frac{c_0}{2}\right)\overline{u}}
$
in $E_{4R}^+\backslash \overline{E}_{R}^+$,
we get  for all $R\geq R_1$,
\[
    u(x)\leq(1-c_0)\left(1+\frac{c_0}{2}\right)\overline{u}
         \leq \left(1-\frac{c_0}{2}\right)\overline{u}
           \quad  \mbox{on }\partial E_{2R}\cap\{x_n\geq0\}.
\]
This implies
\[
    u(x) \leq\left(1-\frac{c_0}{2}\right)\overline{u}
    \quad  \mbox{in}~\overline{\mathbb{R}}^n_+\backslash \overline{E}_{2R_1}^+,
\]
which reaches a contradiction.
\end{proof}

Theorem \ref{Tm-converges} implies the following corollary immediately.

\begin{corollary}\label{Coro-u-convg}
  Let $u$ be as in Theorem \ref{Main-thm}. Then

\begin{equation*}
   u(x)\rightarrow 0\quad \mbox{as }|x|\rightarrow\infty.
\end{equation*}
\end{corollary}
\begin{proof}
  It's obvious and thus we omit it here.
\end{proof}

\subsection{The asymptotic behavior at infinity}\label{S3-2}
In this subsection we obtain the asymptotic behavior at infinity of solutions in Theorem \ref{Main-thm}, through constructing a barrier function.

To get the barrier function, we first give a solution of
\begin{equation}\label{Eq-lap}
  \mathfrak{L}u:=x_n^{2\alpha}\sum_{i=1}^{n-1}D_{ii}u(x)+D_{nn}u(x)=0\quad\mbox{in } \mathbb{R}^{n}_+.
\end{equation}
See \cite{Garofalo-1993-JDE} for more details on the Baouendi-Grushin type operator.
Let
\begin{equation}\label{Def-w}
 w(x',x_n)=\frac{x_n}{\left(|x'|^2+\beta x_n^{2+2\alpha}\right)^\gamma},
\end{equation}
where
 $\beta=\frac{1}{(1+\alpha)^2}$,
 $\gamma=\frac{n-1}{2}+\frac{1}{2(1+\alpha)}$.

Obviously,

\[
 w(x)\thicksim \frac{x_n}{d(x)^Q},
\]
where $d(x)$ and $Q$ are given by \eqref{def-dx}.

Simple calculations deduce that

\begin{align}\label{SZ-Dijw}
  D_iw=&-\frac{2\gamma x_i x_n}{\left(|x'|^2+\beta x_n^{2+2\alpha}\right)^{\gamma+1}}, \quad i<n;\nonumber\\
  D_nw=&\frac{1}{\left(|x'|^2+\beta x_n^{2+2\alpha}\right)^{\gamma}}
        -\frac{\gamma\beta(2+2\alpha) x_n^{2+2\alpha}}
        {\left(|x'|^2+\beta x_n^{2+2\alpha}\right)^{\gamma+1}}; \nonumber\\
  D_{ij}w=&-\frac{2\gamma x_n\delta_{ij}}{\left(|x'|^2+\beta x_n^{2+2\alpha}\right)^{\gamma+1}}
          +\frac{4\gamma(\gamma+1) x_ix_jx_n}{\left(|x'|^2+\beta x_n^{2+2\alpha}\right)^{\gamma+2}},
          \quad i,j<n;\nonumber\\
  D_{in}w=&-\frac{2\gamma x_i }{\left(|x'|^2+\beta x_n^{2+2\alpha}\right)^{\gamma+1}}
          +\frac{2\gamma\beta(2+2\alpha) x_i x_n^{2+2\alpha}}{\left(|x'|^2+\beta x_n^{2+2\alpha}\right)^{\gamma+2}},
       \quad i<n;\nonumber\\
  D_{nn}w=&-\frac{\gamma\beta(2+2\alpha) x_n^{1+2\alpha}}
        {\left(|x'|^2+\beta x_n^{2+2\alpha}\right)^{\gamma+1}}
        -\frac{\gamma\beta(2+2\alpha)^2 x_n^{1+2\alpha}}
        {\left(|x'|^2+\beta x_n^{2+2\alpha}\right)^{\gamma+1}}\\&
        +\frac{\gamma(\gamma+1)\beta^2(2+2\alpha)^2 x_n^{3+4\alpha}}{\left(|x'|^2+\beta x_n^{2+2\alpha}\right)^{\gamma+2}}.\nonumber
\end{align}

Then

\begin{align*}
\mathfrak{L}w
         &=-\frac{2\gamma (n-1)x_n^{1+2\alpha}}{\left(|x'|^2+\beta x_n^{2+2\alpha}\right)^{\gamma+1}}
          +\frac{4\gamma(\gamma+1) |x'|^2x_n^{1+2\alpha}}{\left(|x'|^2+\beta x_n^{2+2\alpha}\right)^{\gamma+2}}
          -\frac{\gamma\beta(2+2\alpha) x_n^{1+2\alpha}}
        {\left(|x'|^2+\beta x_n^{2+2\alpha}\right)^{\gamma+1}}\\
        &\quad-\frac{\gamma\beta(2+2\alpha)^2 x_n^{1+2\alpha}}
        {\left(|x'|^2+\beta x_n^{2+2\alpha}\right)^{\gamma+1}}
        +\frac{\gamma(\gamma+1)\beta^2(2+2\alpha)^2 x_n^{3+4\alpha}}{\left(|x'|^2+\beta x_n^{2+2\alpha}\right)^{\gamma+2}} \\
       &=\frac{\{-2\gamma(n-1)-\gamma\beta(2+2\alpha)(3+2\alpha)\}x_n^{1+2\alpha} }
         {\left(|x'|^2+\beta x_n^{2+2\alpha}\right)^{\gamma+1}}
       +\frac{4\gamma(\gamma+1)\{|x'|^2+\beta^2(1+\alpha)^2 x_n^{2+2\alpha}\}x_n^{1+2\alpha}}
        {\left(|x'|^2+\beta x_n^{2+2\alpha}\right)^{\gamma+2}}\\
      &\xlongequal{\beta=(1+\alpha)^{-2}}
         \frac{\{-2\gamma(n-1)-\gamma\beta(2+2\alpha)(3+2\alpha)\}x_n^{1+2\alpha} }
         {\left(|x'|^2+\beta x_n^{2+2\alpha}\right)^{\gamma+1}}
         +\frac{4\gamma(\gamma+1)x_n^{1+2\alpha}}
        {\left(|x'|^2+\beta x_n^{2+2\alpha}\right)^{\gamma+1}}\\
    &=\frac{2\gamma\{-(n-1)-(1+\alpha)^{-1}(3+2\alpha)+2(\gamma+1)\}x_n^{1+2\alpha} }
        {\left(|x'|^2+\beta x_n^{2+2\alpha}\right)^{\gamma+1}}\\
     &=\frac{2\gamma\{-n+1-(1+\alpha)^{-1}+2\gamma\}x_n^{1+2\alpha} }
        {\left(|x'|^2+\beta x_n^{2+2\alpha}\right)^{\gamma+1}}
    \xlongequal[i.e.,\gamma=\frac{n-1}{2}+\frac{1}{2(1+\alpha)}]{-n+1-(1+\alpha)^{-1}+2\gamma=0}0,
\end{align*}
where $\mathfrak{L}$ is given by \eqref{Eq-lap}.

 {Using} $w$, we can construct a supersolution of \eqref{main-Eq} as follows.

\begin{lemma}\label{supsolution}
Let $L$ be given by
\eqref{main-Eq} with coefficients satisfying \eqref{Sz-unif-ellp}, \eqref{Sz-unif-ellp-2} and \eqref{Sz-aij-good}. Then  for any $\rho\in\left(0,\min\left\{\frac{s}{n-1},1\right\}\right)$,  there exists $R_0\geq1$ depending only on $\rho$, $s$, $\alpha$ and $n$ such that

\begin{equation}\label{Eq-sup-sol-2}
  L(w-w^{1+\rho})\leq0\quad \mbox{in }\mathbb{R}^n_+\backslash \overline{E}_{R_0}^+.
\end{equation}
\end{lemma}

\begin{proof}
For $i,j<n$,

\begin{align}\label{sz-est-Dij}
  |D_{ij}(w^{1+\rho})|=&\left|(1+\rho) w^{\rho} D_{ij}w+\rho(1+\rho) w^{\rho-1}D_iwD_jw\right|\\
  \leq& (1+\rho) w^{\rho}\left\{\frac{2\gamma x_n}{\left(|x'|^2+\beta x_n^{2+2\alpha}\right)^{\gamma+1}}
         +\frac{4\gamma(\gamma+1) |x'|^2x_n}{\left(|x'|^2+\beta x_n^{2+2\alpha}\right)^{\gamma+2}}\right\}
         \nonumber\\
    & +\rho(1+\rho) w^{\rho-1}\frac{4\gamma |x'|^2 x_n^2}
         {\left(|x'|^2+\beta x_n^{2+2\alpha}\right)^{2(\gamma+1)}}
         \nonumber\\
     \leq& \frac{ C(\rho,\alpha,n)w^{\rho}x_n}{\left(|x'|^2+\beta x_n^{2+2\alpha}\right)^{\gamma+1}}
    +\frac{ C(\rho,\alpha,n)w^{\rho-1} x_n^2}
         {\left(|x'|^2+\beta x_n^{2+2\alpha}\right)^{2\gamma+1}}
         \nonumber\\
   \leq&\frac{ C(\rho,\alpha,n)w^{\rho-1}x_n^2}
         {\left(|x'|^2+\beta x_n^{2+2\alpha}\right)^{2\gamma+1}},
         \nonumber
\end{align}
and

\begin{align}\label{sz-est-Din}
 |D_{in}(w^{1+\rho})|=&\left|\rho(1+\rho) w^{\rho-1}D_iwD_nw+(1+\rho) w^{\rho} D_{in}w\right|\\
   \leq&\frac{2\gamma\rho(1+\rho) w^{\rho-1} |x'| x_n}
      {\left(|x'|^2+\beta x_n^{2+2\alpha}\right)^{\gamma+1}}
     \left\{\frac{1}{\left(|x'|^2+\beta x_n^{2+2\alpha}\right)^{\gamma}}
        +\frac{\gamma\beta(2+2\alpha) x_n^{2+2\alpha}}
        {\left(|x'|^2+\beta x_n^{2+2\alpha}\right)^{\gamma+1}}
     \right\}
         \nonumber\\
     &+(1+\rho) w^{\rho}
     \left\{\frac{2\gamma |x'| }{\left(|x'|^2+\beta x_n^{2+2\alpha}\right)^{\gamma+1}}
          +\frac{2\gamma\beta(2+2\alpha) |x'| x_n^{2+2\alpha}}{\left(|x'|^2+\beta x_n^{2+2\alpha}\right)^{\gamma+2}}
          \right\}
         \nonumber\\
     \leq&\frac{ C(\rho,\alpha,n)w^{\rho-1} |x'|x_n}
         {\left(|x'|^2+\beta x_n^{2+2\alpha}\right)^{2\gamma+1}}
         +\frac{ C(\rho,\alpha,n)w^{\rho} |x'|}
         {\left(|x'|^2+\beta x_n^{2+2\alpha}\right)^{\gamma+1}} \nonumber\\
      \leq&\frac{ C(\rho,\alpha,n)w^{\rho-1} |x'|x_n}
         {\left(|x'|^2+\beta x_n^{2+2\alpha}\right)^{2\gamma+1}},
         \nonumber
\end{align}
where $C(\rho,\alpha,n)$ is positive, depending only on $\rho$, $\alpha$ and $n$, and may change from line to line.

Thus,

\begin{align*}
   \mathfrak{L}(w^{1+\rho})
   &=x_n^{2\alpha}\sum_{i=1}^{n-1}(1+\rho) w^{\rho} D_{ii}w+\rho(1+\rho) w^{\rho-1}D_iwD_iw
    +(1+\rho) w^{\rho} D_{nn}w\\
   &\quad+\rho(1+\rho) w^{\rho-1}(D_nw)^2 \\
   &\xlongequal{\mathfrak{L}w=0} x_n^{2\alpha}\rho(1+\rho) w^{\rho-1}\sum_{i=1}^{n-1}(D_iw)^2
    +\rho(1+\rho) w^{\rho-1}(D_nw)^2 \\
   &=\rho(1+\rho) w^{\rho-1}\left\{x_n^{2\alpha}\sum_{i=1}^{n-1}
    \left(
     -\frac{2\gamma x_i x_n}{\left(|x'|^2+\beta x_n^{2+2\alpha}\right)^{\gamma+1}}
    \right)^2\right.\\
   &\left.\quad\quad\quad\quad\quad\quad\quad\quad+
    \left(
       \frac{1}{\left(|x'|^2+\beta x_n^{2+2\alpha}\right)^{\gamma}}
        -\frac{\gamma\beta(2+2\alpha) x_n^{2+2\alpha}}
        {\left(|x'|^2+\beta x_n^{2+2\alpha}\right)^{\gamma+1}}\right)^2\right\}\\
  &=\rho(1+\rho) w^{\rho-1}\left\{\frac{4\gamma^2 |x'|^2 x_n^{2+2\alpha}}
    {\left(|x'|^2+\beta x_n^{2+2\alpha}\right)^{2(\gamma+1)}}
    + \frac{1}{\left(|x'|^2+\beta x_n^{2+2\alpha}\right)^{2\gamma}}
   \right.\\
  &\left.\quad\quad\quad\quad\quad\quad\quad\quad
  -\frac{\gamma\beta(2+2\alpha) x_n^{2+2\alpha}}
        {\left(|x'|^2+\beta x_n^{2+2\alpha}\right)^{2\gamma+1}}
         +\frac{\gamma^2\beta^2(2+2\alpha)^2 x_n^{4+4\alpha}}
        {\left(|x'|^2+\beta x_n^{2+2\alpha}\right)^{2(\gamma+1)}}\right\}\\
  &\xlongequal{\beta=(1+\alpha)^{-2}}
  \rho(1+\rho) w^{\rho-1}
  \left\{
      \frac{4\gamma^2  x_n^{2+2\alpha}}
       {\left(|x'|^2+\beta x_n^{2+2\alpha}\right)^{2\gamma+1}}
       -\frac{2\gamma(1+\alpha)^{-1} x_n^{2+2\alpha}}
        {\left(|x'|^2+\beta x_n^{2+2\alpha}\right)^{2\gamma+1}}\right.\\
     &\left.\quad\quad\quad\quad\quad\quad\quad\quad\quad\quad\quad\quad
       +\frac{1}{\left(|x'|^2+\beta x_n^{2+2\alpha}\right)^{2\gamma}}
  \right\}\\
  &\xlongequal{\gamma=\frac{n-1}{2}+\frac{1}{2(1+\alpha)}}
    \frac{(n-1)\left(n-1+\frac{1}{1+\alpha}\right)\rho(1+\rho) w^{\rho-1} x_n^{2+2\alpha}}
    {\left(|x'|^2+\beta x_n^{2+2\alpha}\right)^{2\gamma+1}}
    +\frac{\rho(1+\rho) w^{\rho-1}}{\left(|x'|^2+\beta x_n^{2+2\alpha}\right)^{2\gamma}},
\end{align*}
where $\mathfrak{L}$ is given by \eqref{Eq-lap}.
This together with \eqref{sz-est-Dij} and \eqref{sz-est-Din} implies

\begin{align}\label{sz-est-wrho-4}
  L \left(w^{1+\rho}\right)
  \geq& \mathfrak{L}\left(w^{1+\rho}\right)
  -\sum_{i,j=1}^{n-1}|a_{ij}(x)-\delta_{ij}||D_{ij}(w^{1+\rho})|x_n^{2\alpha}
     -\sum_{i=1}^{n-1}|a_{in}(x)||D_{in}(w^{1+\rho})|\\
     \geq& \frac{\rho(1+\rho) w^{\rho-1}}{\left(|x'|^2+\beta x_n^{2+2\alpha}\right)^{2\gamma}}
    -\left(|x'|+x_n^{1+\alpha}\right)^{-s}\frac{ C(\rho,\alpha,n)w^{\rho-1}x_n^{2+2\alpha}}
         {\left(|x'|^2+\beta x_n^{2+2\alpha}\right)^{2\gamma+1}}
         \nonumber\\
    &-\left(|x'|+x_n^{1+\alpha}\right)^{-s}\frac{ C(\rho,\alpha,n)w^{\rho-1}|x'|x_n }
         {\left(|x'|^2+\beta x_n^{2+2\alpha}\right)^{2\gamma+1}}
         \nonumber\\
         \geq&\frac{\rho(1+\rho) w^{\rho-1}}{\left(|x'|^2+\beta x_n^{2+2\alpha}\right)^{2\gamma}}
         -\frac{ C(\rho,\alpha,n)w^{\rho-1}}
         {\left(|x'|^2+\beta x_n^{2+2\alpha}\right)^{2\gamma+\frac{s}{2}}}
         -\frac{ C(\rho,\alpha,n)w^{\rho-1}}
         {\left(|x'|^2+\beta x_n^{2+2\alpha}\right)^{2\gamma+\frac{s}{2}+\frac{1}{2}-\frac{1}{2(1+\alpha)}}}
         \nonumber\\
    \geq&\frac{\frac{1}{2}\rho(1+\rho) w^{\rho-1}}
         {\left(|x'|^2+\beta x_n^{2+2\alpha}\right)^{2\gamma}}\quad\mbox{in }\mathbb{R}^n_+\backslash \overline{E}_{R_0}^+
         \nonumber
\end{align}
for some $R_0\geq1$ depending only on $\rho$, $s$, $\alpha$ and $n$.
Similarly,

\begin{align}\label{sz-est-wrho-3}
  L w\leq& \mathfrak{L}w+\sum_{i,j=1}^{n-1}|a_{ij}(x)-\delta_{ij}||D_{ij}w|x_n^{2\alpha}
     +\sum_{i=1}^{n-1}|a_{in}(x)||D_{in}w|\\
     \leq& \frac{ C(\rho,\alpha,n)x_n^{1+2\alpha}}
         {\left(|x'|^2+\beta x_n^{2+2\alpha}\right)^{\gamma+1+\frac{s}{2}}}
    +\frac{ C(\rho,\alpha,n)|x'|}{\left(|x'|^2+\beta x_n^{2+2\alpha}\right)^{\gamma+1+\frac{s}{2}}}
         \nonumber\\
    \leq&\frac{ C(\rho,\alpha,n)x_n^{1+2\alpha}}
         {\left(|x'|^2+\beta x_n^{2+2\alpha}\right)^{\gamma+1+\frac{s}{2}}}
    +\frac{ C(\rho,\alpha,n)|x'|}{\left(|x'|^2+\beta x_n^{2+2\alpha}\right)^{\gamma+1+\frac{s}{2}}}
         \nonumber\\
    \leq&\frac{ C(\rho,\alpha,n)}
         {\left(|x'|^2+\beta x_n^{2+2\alpha}\right)^{\gamma+1+\frac{s}{2}-\frac{1+2\alpha}{2(1+\alpha)}}}.
         \nonumber
\end{align}

Since $\rho\in\left(0,\min\left\{\frac{s}{n-1},1\right\}\right)$, we get

\begin{equation}\label{sz-wrho-exp}
  \frac{w^{\rho-1}}{\left(|x'|^2+\beta x_n^{2+2\alpha}\right)^{2\gamma}}= \frac{x_n^{\rho-1}}{\left(|x'|^2+\beta x_n^{2+2\alpha}\right)^{2\gamma+\gamma(\rho-1)}}
  \geq\left(|x'|^2+\beta x_n^{2+2\alpha}\right)^{-2\gamma-\gamma(\rho-1)-\frac{1-\rho}{2(1+\alpha)}},
\end{equation}
\begin{equation}\label{sz-comp}
  \left(-2\gamma-\gamma(\rho-1)-\frac{1-\rho}{2(1+\alpha)}\right)+ \left(\gamma+1+\frac{s}{2}-\frac{1+2\alpha}{2(1+\alpha)}\right)=-\frac{n-1}{2}\rho+\frac{s}{2}>0.
\end{equation}

By \eqref{sz-est-wrho-4}, \eqref{sz-est-wrho-3}, \eqref{sz-wrho-exp} and \eqref{sz-comp}, we have
\[
L(w-w^{1+\rho})\leq0\quad \mbox{in }\mathbb{R}^n_+\backslash \overline{E}_{R_0}^+
\]
for larger $R_0\geq1$ depending only on $\rho$, $s$, $\alpha$ and $n$.
\end{proof}

\emph{Proof of Theorem \ref{Main-thm}.}

By Lemma \ref{supsolution},
for any fixed $\rho\in\left(0,\min\left\{\frac{s}{n-1},1\right\}\right)$, there exists $R>1$ depending only on  $s$, $\alpha$ and $n$ such that

\begin{equation*}
  L\left(w-w^{1+\rho}\right)\leq0\quad \mbox{in }\mathbb{R}^n_+\backslash \overline{E}_{R}^+.
\end{equation*}

By $u(x)=0$ on $\{x_n=0\}$, $|Du(x)|\leq 1$ in $\mathbb{R}^n_+\backslash \overline{E}_{1}^+$ and Newton-Leibniz formula, we get
\[
|u(x)|\leq 2 x_n\quad \mbox{on }\partial{E}_{R}\cap\{x_n\geq0\}.
\]
On $\partial{E}_{R}\cap\{x_n\geq0\}$, it is clear that

\begin{align*}
   w-w^{1+\rho}=w(1-w^\rho)\geq c(R,\alpha,n)x_n.
\end{align*}
Above two {inequalities} imply that
for some $C>0$ depending only on $s$, $\delta$, $\alpha$ and $n$,

\begin{equation}\label{SZ-u-less-w-1}
    |u(x)|\leq C(w-w^{1+\rho}),
    \quad \mbox{on }\partial{E}_{R}\cap\{x_n\geq0\}.
\end{equation}

For any $\varepsilon>0$, by Corollary \ref{Coro-u-convg},
there exists $R_{\varepsilon}>R$ such that

\begin{equation}\label{SZ-u-less-w-2}
   |u(x)|\leq \varepsilon,\quad x\in\partial E_{R_{\varepsilon}}\cap\{x_n\geq0\}.
\end{equation}
It follows from \eqref{SZ-u-less-w-1}, \eqref{SZ-u-less-w-2}
and $u(x)=0$ on $(E_{R_{\varepsilon}}\backslash \overline{E}_{R})\cap\{x_n=0\}$ that

\begin{equation*}
    |u(x)|\leq C (w-w^{1+\rho})+\varepsilon
    \quad\mbox{on}\;\partial (E^+_{R_{\varepsilon}}\backslash \overline{E}^+_{R}).
\end{equation*}
By the comparison principle,
\[
    |u(x)|\leq C(w-w^{1+\rho})+\varepsilon
    \quad\mbox{in } E^+_{R_{\varepsilon}}\backslash \overline{E}^+_{R}.
\]
Letting $\varepsilon \rightarrow 0$, it follows (\ref{SZ_AsB_Lin}).
$\hfill\Box$

\section{Appendix}

In this section, Lemma \ref{Lm-Unif-Ellp} will be proved as follows.

\begin{proof}
Denote

\[
A'(x)=
  \begin{pmatrix}
  a_{11}(x)& \cdots & a_{1,n-1}(x)\\
  \vdots & \ddots & \vdots \\
  a_{n-1,1}(x) & \cdots &a_{n-1,n-1}(x)
   \end{pmatrix}
\]
and

\[
\widetilde{A}(x)=
  \begin{pmatrix}
   &  & & a_{1,n}(x)x_n^{\alpha}\\
   & A'(x)x_n^{2\alpha} &  &\vdots\\
   &  & &a_{n-1,n}(x)x_n^{\alpha}\\
  a_{n,1}(x)x_n^{\alpha} & \cdots &a_{n,n-1}(x)x_n^{\alpha}&1
   \end{pmatrix},
\]
where $a_{ij}(x)$ and $a_{in}(x)$ are given by \eqref{main-Eq}.

It suffices to show that eigenvalues of $\widetilde{A}(x)$ are positive in $\overline{B}_{1}^+$ and have uniformly bound (depending on the fixed number $\varepsilon_0$) in $\overline{B}_{1}^+\cap\{x_n\geq \varepsilon_0\}$.

If $A'(x)$ have eigenvalues $\lambda_1(x)$,$\cdots$,$\lambda_{n-1}(x)$, by \eqref{Sz-unif-ellp}, we get
$
    \lambda\leq \lambda_i(x)\leq \Lambda$, $i=1,2,\cdots,n-1.
$
It's clear that {there exists a }orthogonal matrix ${P'}_{(n-1)\times(n-1)}$ such that

$$(P')^TA'P'=diag\{\lambda_1(x),\cdots,\lambda_{n-1}(x)\}.$$

Observe that eigenvalues of $\widetilde{A}(x)$ are that of the following matrix
\[
B(x):=P^TAP=
  \begin{pmatrix}
\lambda_1(x)x_n^{2\alpha}&  & & \widetilde{a}_{1,n}(x)x_n^{\alpha}\\
   & \ddots&  &\vdots\\
   &  &\lambda_{n-1}(x)x_n^{2\alpha}&\widetilde{a}_{n-1,n}(x)x_n^{\alpha}\\
  \widetilde{a}_{n,1}(x)x_n^{\alpha} & \cdots &\widetilde{a}_{n,n-1}(x)x_n^{\alpha}&1
   \end{pmatrix}
\]
with
\[
\widetilde{a}_{i,n}(x)=\sum_{j=1}^{n-1}P'_{ij}{a}_{j,n}(x),\quad i=1,\cdots,n-1;   \quad
   P=
   \begin{pmatrix}
  P' & 0\\
  0 & 1
   \end{pmatrix}.
\]

Thus, we only need to show that all eigenvalues of $B(x)$  are positive in $\overline{B}_{1}^+$ and have uniformly  bound in $\overline{B}_{1}^+\cap\{x_n\geq \varepsilon_0\}$.
Let
\[
e_i=(\underbrace{0,\cdots,0,1}_{i},0,\cdots,0),\quad i\leq n.
\]
we have

\begin{equation*}\label{Sz-eiei}
   e_i^TB(x)e_i=\lambda_ix_n^{2\alpha}, \quad
   e_i^TB(x)e_n=a_{in}x_n^{\alpha},\quad i\leq n-1;\quad
  e_n^TB(x)e_n=1.
\end{equation*}

For any $\xi\in\mathbb{R}^n$ with $|\xi|=1$, there exists a unique consequence $\{b_i\}_{i=1}^{n}$ such that
\[
\xi=\sum_{i=1}^{n}b_ie_i,\quad \sum_{i=1}^{n} b_i^2=1.
\]
Then

\begin{align*}
   \xi^TB(x)\xi&=\sum_{i=1}^{n}(b_ie_i)^TB(x)(b_ie_i)
   =\sum_{i=1}^{n}b_i^2e_i^TB(x)e_i\\
   &=\sum_{i=1}^{n-1}\lambda_i b_i^2 x_n^{2\alpha}
      +\sum_{i=1}^{n-1}2b_i b_n\widetilde{a}_{i,n}x_n^{\alpha}+b_n^2\\
   &\geq \lambda x_n^{2\alpha}\sum_{i=1}^{n-1}b_i^2
         +\sum_{i=1}^{n-1}2 b_i b_n \widetilde{a}_{i,n}x_n^{\alpha} +b_n^2
        \quad(\mbox{by }\eqref{Sz-unif-ellp}).
\end{align*}
Applying Cauchy's inequality to $2 b_i b_n \widetilde{a}_{i,n}x_n^{\alpha}$, we have that for any $\tau\in(0,1)$,

\begin{align*}
      \left|\sum_{i=1}^{n-1}2 b_i b_n \widetilde{a}_{i,n}x_n^{\alpha}\right|&\leq \tau\sum_{i=1}^{n-1}
      \left\{\lambda^{\frac{1}{2}}b_i x_n^{\alpha}\right\}^2
      +\tau^{-1}\sum_{i=1}^{n-1}\left\{\lambda^{-\frac{1}{2}}b_n \widetilde{a}_{i,n}\right\}^2\\
      &=\tau\lambda x_n^{2\alpha}\sum_{i=1}^{n-1}b_i^2
      +\tau^{-1} b_n^2\lambda^{-1}\sum_{i=1}^{n-1}\widetilde{a}_{i,n}^2.
  \end{align*}
Therefore, for any $\tau\in(1-\delta,1)$,

\begin{align*}
     \xi^TB(x)\xi&\geq \lambda x_n^{2\alpha}\sum_{i=1}^{n-1}b_i^2+b_n^2
     -\tau\lambda x_n^{2\alpha}\sum_{i=1}^{n-1}b_i^2
      -\tau^{-1}b_n^2\lambda^{-1}\sum_{i=1}^{n-1}\widetilde{a}_{i,n}^2\\
        &\geq (1-\tau)\lambda x_n^{2\alpha}\sum_{i=1}^{n-1}b_i^2
        +b_n^2\left\{1-\tau^{-1}(1-\delta)\right\}\qquad (\mbox{by }\eqref{Sz-unif-ellp-2}),
\end{align*}
which implies {that} $L$ is elliptic in $\overline{B}_{1}^+$. And then if $\{x_n\geq \varepsilon_0\}$,

\begin{align*}
        \xi^TA(x)\xi&\geq(1-\tau)\lambda \varepsilon_0^{2\alpha}\sum_{i=1}^{n-1}b_i^2
        +b_n^2\{1-\tau^{-1}(1-\delta)\}
        \geq\min\left\{(1-\tau)\lambda \varepsilon_0^{2\alpha},1-\tau^{-1}(1-\delta)\right\}.
  \end{align*}
In particular, taking $\tau=1-\frac{1}{2}\delta$, we have that for any  $x\in \overline{B}_{1}^+\cap\{x_n\geq\varepsilon_0\}$,
\[
    \xi^TA(x)\xi\geq \min\left\{\frac{1}{2}\delta\lambda \varepsilon_0^{2\alpha},1-\left(1-\frac{1}{2}\delta\right)^{-1}(1-\delta)\right\}>0.
\]
Therefore, eigenvalues of $B(x)$ have uniformly below bound in $\overline{B}_{1}^+\cap\{x_n\geq \varepsilon_0\}$.

Similarly, one can also obtain the uniformly upper bound of eigenvalues of $B(x)$ in $\overline{B}_{1}^+\cap\{x_n\geq \varepsilon_0\}$.
\end{proof}

%

\end{document}